\documentclass[12pt,reqno]{amsart}
\usepackage[activeacute,english]{babel}

\usepackage[utf8]{inputenc}

\usepackage[english]{babel}
\usepackage{amsfonts,amsmath,amsthm,bm,amssymb}
\usepackage{graphics,tikz,float}

\usepackage{hyperref,cleveref}
\usepackage{autonum}

\DeclareRobustCommand{\SkipTocEntry}[5]{}

\usepackage[margin=0.9in]{geometry}
\newcommand{\R}{{\mathbb R}}
\newcommand{\Z}{{\mathbb Z}}

\newcommand{\sgn}{{\operatorname{sign}}}

\newtheorem{theorem}{Theorem}[section]
\newtheorem{lemma}[theorem]{Lemma}

\newtheorem{proposition}[theorem]{Proposition}

\numberwithin{equation}{section}

\allowdisplaybreaks[1]


\begin{document}

\title[Optimal H\"older regularity]
{Examples of optimal  H\"older regularity in semilinear equations involving the fractional Laplacian}

\author[G. Csat\'o]{Gyula Csat\'o}
\address{G. Csat\'o \textsuperscript{1,2}
\newline
\textsuperscript{1}
Departament de Matem\`{a}tiques i Inform\`{a}tica,
Universitat de Barcelona,
Gran Via 585,
08007 Barcelona, Spain
\newline
\textsuperscript{2}
Centre de Recerca Matem\`atica, Edifici C, Campus Bellaterra, 08193 Bellaterra, Spain}
\email{gyula.csato@ub.edu}

\author[A. Mas]{Albert Mas}
\address{A. Mas \textsuperscript{1,2}
\newline
\textsuperscript{1}
Departament de Matem\`atiques,
Universitat Polit\`ecnica de Catalunya,
Campus Diagonal Bes\`os, Edifici A (EEBE), Av. Eduard Maristany 16, 08019
Barcelona, Spain
\newline
\textsuperscript{2}
Centre de Recerca Matem\`atica, Edifici C, Campus Bellaterra, 08193 Bellaterra, Spain}
\email{albert.mas.blesa@upc.edu}

\date{\today}
\thanks{The two authors are supported by the Spanish grants PID2021-123903NB-I00 and RED2022-134784-T funded by MCIN/AEI/10.13039/501100011033 and by ERDF ``A way of making Europe'', and by the Catalan grant 
2021-SGR-00087. The first author is in addition supported by the Spanish grant PID2021-125021NA-I00.
This work is supported by the Spanish State Research Agency, through the Severo Ochoa and Mar\'ia de Maeztu Program for Centers and Units of
Excellence in R\&D (CEX2020-001084-M)}

\begin{abstract}
We discuss the H\"older regularity of solutions to  the semilinear equation involving the fractional Laplacian $(-\Delta)^s u=f(u)$ in one dimension. 
We put in evidence a new regularity phenomenon which is a  combined effect of the nonlocality and the semilinearity of the equation, since it does not happen neither for local semilinear equations, nor for nonlocal linear equations. 
Namely, for nonlinearities $f$ in $C^\beta$ and when $2s+\beta <1$, the solution is not always $C^{2s+\beta-\epsilon}$ for all $\epsilon >0$. Instead, in general the solution $u$ is at most $C^{2s/(1-\beta)}.$

\bigskip

\end{abstract}

\maketitle



\section{Introduction}

It is well known that if $f\in C^{\beta}(\R)$ and $\beta>0$ is not an integer, then bounded solutions to $-\Delta u=f(u)$ in $B_1\subset\R^n$ satisfy  $u\in C^{2+\beta}(B_{1/2})$, a result which is obtained by simply combining Calder\'on-Zygmund and Schauder's estimates.
It can be easily seen that this regularity result is optimal. In  \Cref{section:appendix opt reg local} we give a simple example for this fact in dimension $1$, since it is the starting point for one of our examples in the nonlocal case.
If $\beta\geq0$ is an integer and $n>1$, then $u$ is of class $C^{2+\beta-\epsilon}$ for all $\epsilon>0$, but to prove or disprove that $u$ is always of class $C^{2+\beta}$ still seems to be an unsettled problem. In the case $\beta=0,$ and thus $f$ is merely continuous, it has been posed as an open question in \cite{Shahgholian 2015} whether all solutions have continuous second order derivatives.

In this paper we discuss the H\"older regularity for semilinear equations involving the fractional Laplacian in dimension 1, that is, the regularity of solutions $u$ to the equation  
\begin{equation}
 \label{eq:frac Lapl intro}
  (-\Delta)^s u=f(u) \quad\text{ in }(-1,1),
\end{equation}
where $f$ is a given H\"older function, $0<s<1$, and 
\begin{equation}\label{defi laplacian_}
(-\Delta)^s u(x):=c_s\, \text{P.V.}\!\int_\R dy\,
\frac{u(x)-u(y)}{|x-y|^{1+2s}},\qquad
c_s:=\frac{s4^s\Gamma(1/2+s)}{\sqrt\pi\,\Gamma(1-s)}.
\end{equation}
Recall that if $u\in C^{\alpha}(\R)$
for some $\alpha>2s$ then $(-\Delta)^s u(x)$ is well defined and finite for every $x\in (-1,1)$.\footnote{Given 
$\Omega\subset\R$ and $\alpha\geq0$, throughout this paper we denote by $C^\alpha(\Omega)$ the space of functions defined on $\Omega$ that are continuously differentiable  with bounded derivatives up to order $\lfloor\alpha\rfloor$ (continuous and bounded for $\lfloor\alpha\rfloor=0$) and whose 
$\lfloor\alpha\rfloor$-th derivatives satisfy the H\"older condition with index $\alpha-\lfloor\alpha\rfloor$ 
(we denote by $\lfloor\alpha\rfloor$ the integer part 
of~$\alpha$). In particular, we always have $C^\alpha(\Omega)\subset L^\infty(\Omega)$.}

In view of the local result, one would expect that bounded solutions to $(-\Delta)^s u=f(u)$ with $f\in C^{\beta}(\R)$ for some $\beta>0$ are of class $C^{2s+\beta-\epsilon}$ for every $\epsilon>0$ (or even perhaps $C^{2s+\beta}$), since the order of the operator $(-\Delta)^s$ is $2s$. 
In our recent paper \cite{CCM Semilinear Var} we have proven that this claim is indeed true when $2s\geq1-\beta$. However,  if $2s<1-\beta$ we were only able to prove that bounded solutions are of class $C^{2s/(1-\beta)}$. The exponent $2s/(1-\beta)$ is strictly smaller than $2s+\beta$ when $2s<1-\beta$, and seems to be new in the literature on fractional semilinear elliptic equations.
Let us also mention that in the linear case 
$(-\Delta)^su=h$ with 
$h\in C^\beta(\R)$ (i.e., the right hand side of the equation is independent of $u$), the solution $u$ is of class $C^{2s+\beta}$ and the H\"older exponent $2s+\beta$ is known to be sharp for the solution, as has been shown in \cite[Corollary 5.4]{Bass}. This fact in the linear case is independent of whether $2s$ is greater or smaller than $1-\beta.$

We now cite the result obtained in \cite{CCM Semilinear Var} where, for the first time, this new H\"older exponent $2s/(1-\beta)$ was discovered. The following version is a simplification of the statements in \cite{CCM Semilinear Var}, which deal with more general integro-differential operators and distributional solutions, in both a periodic and nonperiodic setting.

\begin{theorem}$($\cite[Theorem 1.7 and Remark 8.2]{CCM Semilinear Var}$)$\label{fthm:Holder reg periodic} 
Let $0<s<1$, $f\in C^\beta(\R)$ for some $\beta>0$, and 
$u\in C^{\alpha}(\R)$ for some $\alpha>2s$. Assume that
$$
(-\Delta)^s u=f(u) \quad\text{ in }(-1,1).
$$
Then, the following holds:
\begin{itemize}
\item[$(i)$] If $2s\leq1-\beta$ then $u\in C^{\frac{2s}{1-\beta}-\epsilon}([-\frac1 2,\frac1 2])$ for all $\epsilon>0$.

\item[$(ii)$] If $2s>1-\beta$ then $u\in C^{\beta+2s-\epsilon}([-\frac1 2,\frac1 2])$ for all $\epsilon>0$. 
\end{itemize}
\end{theorem}

To see that $(i)$ is a natural result, consider the simple heuristic calculation: assume that  $r<1$ and
$(-\Delta)^s(|x|^r)\approx |x|^{r-2s}=(|x|^r)^{(r-2s)/r}=:f(|x|^r)$. Here we have $\beta =1-2s/r$ and, thus, 
$r=2s/(1-\beta)$.  Besides this heuristic argument, the precise reason why the corresponding proof of $(ii)$ (or its analogue in the local case) breaks down
in case $(i)$ 
is that the composition $f(v)$ of H\"older functions $f\in C^{\beta}$ and $v\in C^{\sigma}$ behaves differently when $\sigma$ and $\beta$ are both smaller than~$1$, and one merely obtains $f(v)\in C^{\beta\sigma}.$ In the iteration argument of the proof of \Cref{fthm:Holder reg periodic}~$(i)$ one uses in each step this composition argument, never reaching a H\"older exponent greater than or equal to $1$.

The~$\epsilon$ appearing in the statement of \Cref{fthm:Holder reg periodic} $(ii)$ is due to the fact that $\beta+2s$ may be an integer. In particular, 
if $\beta+2s$ is not an integer, and if in the iteration argument of the proof  of \Cref{fthm:Holder reg periodic}~$(ii)$ no H\"older exponent ever is an integer, then we indeed get $u\in C^{\beta+2s}([-\frac1 2,\frac1 2])$. 
On the contrary, the~$\epsilon$ appearing in the statement of \Cref{fthm:Holder reg periodic}~$(i)$ has a different nature. It comes from the fact that the iteration method of the proof yields 
$u$ of class $C^{2s\beta_{k}}$ for all $k\geq 1$, where 
$\beta_k$ is a sequence of exponents such that 
$\beta_k\uparrow1/(1-\beta)$ as $k\uparrow+\infty$ but 
$\beta_k<1/(1-\beta)$ for all $k$. 

The main goal of the present paper is to show that in the case of \Cref{fthm:Holder reg periodic}~$(i)$ one indeed cannot go beyond $C^{2s/(1-\beta)}$ regularity. Currently we do not know whether \Cref{fthm:Holder reg periodic}~$(i)$ is sharp or if, instead, it also holds true for 
$\epsilon=0$ as well; we leave this question as an open problem. But what we know for sure is that the exponent $2s/(1-\beta)$ cannot be improved, as shown in the next theorem. For the sake of completeness, we also deal with the optimality of the exponent in the case of \Cref{fthm:Holder reg periodic}~$(ii)$. Our result refers to the periodic setting since it was announced in \cite[Section 1.4]{CCM Semilinear Var} in the context of periodic solutions but, obviously, it also shows the optimality of the exponents in the Dirichlet case described in \Cref{fthm:Holder reg periodic}.

\begin{theorem}\label{gy:thm:intro:examples reg}
Let $0<s<1$, $0<\beta<1$, and $L>0$. Consider the equation
\begin{equation}\label{gy:intro:Deltas u is fu}
   (-\Delta)^su=f(u)\quad\text{ in }\R.
\end{equation}
\begin{itemize}
\item[$(i)$] If $2s\leq1-\beta$ then there exists 
 $u\in C^{\frac{2s}{1-\beta}}(\R)$
 which is $2L$-periodic, satisfies \eqref{gy:intro:Deltas u is fu} for some $f\in C^{\beta}(\R),$ and
$u\notin C^{\frac{2s}{1-\beta}+\epsilon}([-\rho,\rho])$ for any $\epsilon,\rho>0$.

\item[$(ii)$] If $2s>1-\beta$ then  there exists $u\in C^{2s+\beta}(\R)$ which is $2L$-periodic, satisfies \eqref{gy:intro:Deltas u is fu} for some $f\in C^{\beta}(\R),$ and
$u\notin C^{2s+\beta+\epsilon}([-\rho,\rho])$ for any $\epsilon,\rho>0$.

\end{itemize}
\end{theorem}

The main idea behind this result has nothing to do with periodicity. This is why the proof of \Cref{gy:thm:intro:examples reg} begins by constructing an example where \eqref{gy:intro:Deltas u is fu} holds only in a bounded interval instead of $\R$, and with bounded Dirichlet data on the complement of the interval instead of periodic boundary conditions. This case is simpler and will make the exposition more comprehensible. 

Regarding the structure of the paper, in  \Cref{section:nonperiodic counterexample i}  we deal with the nonperiodic case and $2s\leq 1-\beta$. Then, in  \Cref{section:periodic counterexample i} we modify the arguments to deal with the periodic setting, giving the proof of \Cref{gy:thm:intro:examples reg} $(i)$ ---the much simpler case $2s=1-\beta$ is dealt with separately; see \Cref{lemma:regularity example no surprise} for the nonperiodic setting and the comment below \Cref{gy:example_8periodic:lemma} for the periodic case.
Next, in  \Cref{section:nonperiodic counterexample ii} we address the case $2s>1-\beta$, both in the nonperiodic and the periodic settings, giving the proof of \Cref{gy:thm:intro:examples reg} $(ii)$. Finally, in  \Cref{section:appendix opt reg local} give a simple example for the local semilinear case. In all these situations we will denote the functions at hand by $u$,  since there is no risk of confusion, but they refer to different functions in the different situations.

\section{The nonperiodic example for $2s\leq 1-\beta$}
\label{section:nonperiodic counterexample i}
Let us define the function $u$ by 
\begin{equation}\label{eq:gy:u def}
  u(x):=\left\{\begin{array}{rl}
        -1 & \text{if }x< -1, \\
        -|x|^r & \text{if } -1\leq x< 0,
        \\
        x^r & \text{if }0\leq x\leq 1,
        \\
        1 & \text{if }1< x,
        \end{array}\right.
\end{equation}
where $r>0$.
In this section we will show that, for $r=2s/(1-\beta)$, the function $u$ given by \eqref{eq:gy:u def} is a suitable example for the nonperiodic setting and $2s< 1-\beta$; this is the content of the following result. The case $2s=1-\beta$ is treated separately in \Cref{lemma:regularity example no surprise}.

\begin{proposition}
\label{prop:gy:counterexample reg}
Let $0<s<1$ and $0<\beta<1$ such that $2s< 1-\beta$, and let $u$ be as in \eqref{eq:gy:u def} with $r=2s/(1-\beta)$.
Then,
 there exists 
 $f\in C^{\beta}(\R)$ such that 
$$
  (-\Delta)^su(x)=f(u(x))\quad\text{for all }x\in 
  {\textstyle [-\frac{1}{2},\frac{1}{2}]}.
$$
Moreover, $u$ satisfies
$$
  u\in C^{\frac{2s}{1-\beta}}(\R)\quad\text{and}\quad u\notin C^{\frac{2s}{1-\beta}+\epsilon}([-\rho,\rho])
$$
for any $\epsilon,\rho>0$. 
\end{proposition}

The last two statements concerning the regularity of $u$ are obvious. The difficulty is to prove the existence of $f$ with the required regularity and that $(-\Delta)^su=f(u)$ in $\textstyle [-\frac{1}{2},\frac{1}{2}]$.

Let us start with some preliminary results which do not require this particular choice for $r$ in \Cref{prop:gy:counterexample reg}. The following lemma will be used in the study of the fractional Laplacian of $u$.

\begin{lemma}
\label{lemma:gy:technical lemma for H Lip}
Assume that $0<r<1$ and 
$
  0\leq x\leq z\leq 2^{-r}.
$
Then,
\begin{itemize}
\item[$(i)$]
     $(z^{1/r}+y)^r-(x^{1/r}+y)^r\leq 2^r(z-x)$\quad if $y\geq z^{1/r}$,
\item[$(ii)$]
      $(y-x^{1/r})^r-(y-z^{1/r})^r
      \leq 2^{r-1}(y-z^{1/r})^{r-1}(z-x)$
      \quad if $y>z^{1/r}$.
\end{itemize}
In particular, $(i)$ holds for all $y\geq\frac{1}{2}.$
\end{lemma}

\begin{proof}
To prove $(i)$, consider the function $g$ defined by
$$
   g(y):= (z^{1/r}+y)^r-(x^{1/r}+y)^r-2^r(z-x).
$$
Since $x\leq z$, we have
$
g(z^{1/r})= -(x^{1/r}+z^{1/r})^r+2^rx\leq0.
$
Moreover, using again that $x\leq z$, and also that $0<r<1$, we get
$$
   g'(y)= r\big((z^{1/r}+y)^{r-1}-(x^{1/r}+y)^{r-1}\big)\leq0.
$$
Therefore, $g(y)\leq0$ for all $y\geq z^{1/r}$, which proves $(i)$.

For the proof of $(ii)$ we can assume that $x<z$. If we divide the inequality by $z-x$, the left hand side can be written as
$$
  \frac{(y-x^{1/r})^r-(y-z^{1/r})^r}
  {z-x}=
  \frac{(y-x^{1/r})^r-(y-z^{1/r})^r}{z^{1/r}
  -x^{1/r}}\cdot \frac{z^{1/r}-x^{1/r}}{z-x}=:
  A\cdot B,
$$
where $A$ and $B$ abbreviate the two quotients on the right hand side. By the mean value theorem we have that $B\leq \frac{1}{r}x_0^{(1-r)/r}$ for some $x_0\in (0,2^{-r}),$ and hence
$$
  B\leq \frac{2^{r-1}}{r}.
$$ 
To bound $A$, we apply the mean value theorem to the function  $\psi(t):=(y-t)^r$ and obtain that, for some 
$t_0\in (x^{1/r},z^{1/r})$,
$$
  A=-\frac{(y-z^{1/r})^r-(y-x^{1/r})^r}{z^{1/r}
  -x^{1/r}}=-\psi'(t_0)
  =r(y-t_0)^{r-1}\leq r(y-z^{1/r})^{r-1}.
$$
This concludes the proof of $(ii)$.
\end{proof}

The next lemma studies the regularity properties of  
$(-\Delta)^su$ near the origin. This is the main ingredient in the proof of \Cref{prop:gy:counterexample reg}. In view of its application in other sections, we state the first part of the lemma for more general functions instead of $u$ defined in \eqref{eq:gy:u def}.

\begin{lemma}
\label{lemma:gy:fracLaplace of u}
Given $0<s<1$ and $0<r<1+2s$, let $v:\R\to\R$ be a bounded measurable function such that 
$$
    v(x)=\left\{\begin{array}{rl}
        -|x|^r & \text{if } -1\leq x< 0,
        \\
        x^r & \text{if } 0\leq x\leq 1.
        \end{array}\right.
$$
Then, there exist two constants $c_1,c_2\in\R$ depending only on $r$ and $s$ such that
\begin{equation}
\label{eq:gy:fracLaplace of u}
 (-\Delta)^sv(x)=c_1 x^{r-2s}+c_2 x^r+c_sG(x)-c_sH(x)\quad
 \text{ for all }x\in {\textstyle (0,\frac{1}{2}]},
\end{equation}
where $c_s$ is the constant given in the definition of 
$(-\Delta)^s$ in \eqref{defi laplacian_}, 
\begin{equation}
\begin{split}
\label{eq:gy:G def}
  G(x)&:=x^{r-2s}\int_{\frac{1}{2x}}^{+\infty}dt\,
  \frac{(t+1)^r-(t-1)^r}{t^{1+2s}}\quad\text{for }x>0,
  \\
  \label{eq:gy:H def}
  H(x)&:=\int_{\R\setminus(-\frac{1}{2},\frac{1}{2})}
  dy\,\frac{v(x+y)-v(x)}{|y|^{1+2s}}\quad\text{for }x\geq0.
\end{split}
\end{equation}
Moreover, setting $G(0):=0$ it holds that
\begin{align}
\label{eq:gy:G property}
  G\in C^1({\textstyle [0,\frac{1}{2}]}),
\end{align}
and if $r<1$ and $v=u$ is as in \eqref{eq:gy:u def} then 
\begin{align}
  \label{eq:gy:H property}
  x\mapsto H(x^{1/r})\quad\text{is a Lipschitz function in }[0,2^{-r}].
\end{align}
\end{lemma}

A remark concerning  \eqref{eq:gy:H property} is in order.  It is immediate to check that $H\in C^r({\textstyle [0,\frac{1}{2}]})$, but solely this does not yield that $H\circ\psi$ is Lipschitz in $[0,2^{-r}]$ for every $\psi\in C^{1/r}([0,2^{-r}])$ with values in ${\textstyle [0,\frac{1}{2}]}$. However, since $H$ behaves like $x\mapsto x^{r}$ for $x$ close to 0, if we choose 
$\psi(x)=x^{1/r}$ then $H\circ\psi$ turns out to behave like a linear function, and hence being locally Lipschitz, as \eqref{eq:gy:H property} claims.

\begin{proof}[Proof of \Cref{lemma:gy:fracLaplace of u}]
We first note that the assumption $r<1+2s$ guarantees that the  integrand defining $G$ is in $L^1(a,+\infty)$ for every $a>0.$ Indeed, using the Taylor expansion
$(1+z)^r=1+rz+O(z^2)$ as $z\to 0$ we get
\begin{equation}\label{gy:eq:Taylor exp r}
  (t+1)^r-(t-1)^r=t^r\left(\Big(1+\frac{1}{t}\Big)^r-
  \Big(1-\frac{1}{t}\Big)^r\right)
  =t^{r-1}\big(2r+O(t^{-1})\big)
\end{equation}
as $t\uparrow+\infty$, which proves that $G$ is well defined. This also justifies the upcoming computations dealing with the functions $F$ and $\Phi$ introduced below.

Let us now start with the proof. We clearly have that
$$
  (-\Delta)^sv(x)=c_s\, \text{P.V.}\!\int_{-1/2}^{1/2}dy\, 
  \frac{v(x)-v(x+y)}{|y|^{1+2s}}-c_sH(x)=:c_sA(x)-c_sH(x),
$$
where $A(x)$ is simply an abbreviation for the first integral over ${\textstyle (-\frac{1}{2},\frac{1}{2})}$. Now, if $x\in{\textstyle (0,\frac{1}{2}]}$ and $y\in {\textstyle (-\frac{1}{2},\frac{1}{2})}$, then $x+y\in (-1,1).$ Setting $y=x t$ gives, by the definition of $v$, that
$$
  v(x+y)=v(x(1+t))=x^r\varphi(1+t),
  \quad\text{ where }
  \varphi(z):=\left\{\begin{array}{rl}
         -|z|^r & \text{if }z<0,
         \\
         z^r & \text{if }z\geq0.
         \end{array}\right.
$$
By the change of variable $y=xt$, we get 
\begin{equation}
\begin{split}
  -A(x)&=x^{-2s}\text{ P.V.}\!\int_{-\frac{1}{2x}}^{\frac{1}{2x}}dt\,  
  \frac{v(x(1+t))-v(x)}{|t|^{1+2s}}\\
  &=
  x^{r-2s}\text{ P.V.}\!\int_{-\frac{1}{2x}}^{\frac{1}{2x}}dt\, 
  \frac{\varphi(1+t)-1}{|t|^{1+2s}}
  =x^{r-2s}F\Big(\frac{1}{2x}\Big),
\end{split}
\end{equation}
where
$$
  F(a):=\text{P.V.}\!\int_{-a}^a dt\, \frac{\varphi(1+t)-1}{|t|^{1+2s}}
  \quad\text{for $a\geq1$ (recall that ${\textstyle 0<x\leq\frac{1}{2}}$)}.
$$
We now write $F(a)=Q_1(a)+Q_2+Q_3(a),$ where
$$
  Q_1(a):=\!\int_{-a}^{-1} dt\, \frac{\varphi(1+t)-1}{|t|^{1+2s}},
  \quad\!
  Q_2:=\text{P.V.}\!\int_{-1}^{1} dt\, \frac{\varphi(1+t)-1}{|t|^{1+2s}},
  \quad\!
  Q_3(a):=\!\int_{1}^{a} dt\, \frac{\varphi(1+t)-1}{|t|^{1+2s}}.
$$
Using that the integrand defining $Q_2$ is equal to $((1+t)^r-1)|t|^{-1-2s}$ we easily see that the principal value is well defined and that $Q_2\in \R$ (moreover, if $2s<1$ then there is no need of $\text{P.V.}$ in the definition of $Q_2$).  
Regarding $Q_1(a)$, since $\varphi(1+t)=-|1+t|^r$ if $1+t<0,$ by the change of variable $t\mapsto -t$ we obtain that
$$
  Q_1(a)=-\int_1^a dt\,\frac{(t-1)^r+1}{t^{1+2s}}.
$$
Therefore,
\begin{equation}
\begin{split}
  F(a)&=Q_2+\int_{1}^{a} dt\, 
  \frac{(t+1)^r-(t-1)^r-2}{t^{1+2s}}\\
  &= Q_2-2\int_1^a \frac{dt}{t^{1+2s}}
  +
  \int_{1}^{a} dt\, 
  \frac{(t+1)^r-(t-1)^r}{t^{1+2s}}
  \\
  &= Q_2+\frac{1}{s}(a^{-2s}-1)
  +
  \int_{1}^{+\infty} dt\, 
  \frac{(t+1)^r-(t-1)^r}{t^{1+2s}} 
  -\int_{a}^{+\infty} dt\, 
  \frac{(t+1)^r-(t-1)^r}{t^{1+2s}} 
  \\ 
  &= 
  \overline{Q_2}+ \frac{a^{-2s}}{s}-\int_{a}^{+\infty} dt\, 
  \frac{(t+1)^r-(t-1)^r}{t^{1+2s}}
\end{split}
\end{equation}
for some constant $\overline{Q_2}\in\R$; in the last equality above we used \eqref{gy:eq:Taylor exp r} and that $r<1+2s$.
Thus, we finally obtain that 
$$
  A(x)=-x^{r-2s}F\Big(\frac{1}{2x}\Big)
  =
  -\overline{Q_2}x^{r-2s}-\frac{2^{2s}}{s}x^{r}
  +x^{r-2s}\int_{\frac{1}{2x}}^{+\infty} dt\, 
  \frac{(t+1)^r-(t-1)^r}{t^{1+2s}},
$$
which proves \eqref{eq:gy:fracLaplace of u} for all $0<x\leq{\textstyle \frac{1}{2}}$.

We will now prove \eqref{eq:gy:G property}. Using \eqref{gy:eq:Taylor exp r} we clearly see that $G\in C^1((0,1/2])\cap C([0,1/2])$. Thus, to get 
\eqref{eq:gy:G property} we only have to show that $\lim_{x\downarrow 0}G'(x)$ exists and is finite. We shall use the notation
\begin{equation}
  G(x)=x^{r-2s} \Phi\Big(\frac{1}{2x}\Big), \quad\text{where}\quad 
  \Phi(a):=\int_a^{+\infty}dt\, \frac{(t+1)^r-(t-1)^r}{t^{1+2s}}
  \quad\text{for }a\geq 1,
\end{equation}
and
\begin{equation}
  G'(x)=-\frac{1}{2}\,x^{r-2s-2}\Phi'\Big(\frac{1}{2x}\Big)
  +
  (r-2s)x^{r-2s-1}\Phi\Big(\frac{1}{2x}\Big)=:\Sigma_1(x)+\Sigma_2(x).
\end{equation}
Setting $a=1/(2x),$ we obtain
\begin{align}
  \lim_{x\downarrow 0}\Sigma_1(x)
  =
  c_3\lim_{a\uparrow+\infty} a^{-r+2s+2}\frac{(a+1)^r-(a-1)^r}{a^{1+2s}}
  =c_3\lim_{a\uparrow+\infty} a^{1-r}\left((a+1)^r-(a-1)^r\right)
\end{align}
for some constant $c_3$ depending only on $r$ and $s.$ 
Using now the Taylor expansion \eqref{gy:eq:Taylor exp r} proves that 
$\lim_{x\downarrow 0}\Sigma_1(x)=2rc_3$ exists and is finite.

To compute the limit of $\Sigma_2(x)$ as $x\downarrow 0$ we use the l'H\^opital rule to obtain that, for some constant $c_4$ depending only on $r$ and $s$ that may be renamed in some steps if necessary,
\begin{align}
  \lim_{x\downarrow 0}\Sigma_2(x)
  =
  c_4\lim_{x\downarrow 0}\frac{\Phi\big(\frac{1}{2x}\big)}{x^{1-r+2s}}
  =
  c_4\lim_{x\downarrow 0}\Phi'\Big(\frac{1}{2x}\Big)x^{r-2s-2}
  =c_4\lim_{x\downarrow 0}\Sigma_1(x).
\end{align}
Thus, $\lim_{x\downarrow 0}\Sigma_2(x)$ exists too, and this completes the proof of \eqref{eq:gy:G property}.

We will now prove \eqref{eq:gy:H property}. We decompose $H$ as $H(x)=H_+(x)+H_-(x),$ where
\begin{equation}\label{eq:gy:decomp of H}
  H_+(x):=\int_{\frac{1}{2}}^{+\infty}dy\, \frac{u(x+y)-u(x)}{y^{1+2s}}
  \quad\text{ and }\quad
  H_-(x):=\int_{-\infty}^{-\frac{1}{2}}dy\, \frac{u(x+y)-u(x)}{|y|^{1+2s}}.
\end{equation}
Let us first show that $x\mapsto H_+(x^{1/r})$ is Lipschitz in $[0,2^{-r}].$ Assume that
\begin{equation}\label{eq:gy:x smaller z}
  0\leq x<z\leq 2^{-r},\quad\text{and thus }\quad
  0\leq x^{{1}/{r}}<z^{{1}/{r}}\leq\frac{1}{2}.
\end{equation}
Then,
\begin{equation}\label{eq:gy:Hplus Lipschitz 0}
  \left|H_+(x^{{1}/{r}})-
  H_+(z^{{1}/{r}})\right|\leq
  \int_{\frac{1}{2}}^{+\infty}dy\,
  \frac{\big|u(x^{{1}/{r}}+y)-u(z^{{1}/{r}}+y)\big|}
  {y^{1+2s}}
  +
  \int_{\frac{1}{2}}^{+\infty}dy\,
  \frac{\big|u(x^{{1}/{r}})-u(z^{{1}/{r}})\big|}
  {y^{1+2s}}.
\end{equation}
By the definition of $u$ in  \eqref{eq:gy:u def} we have that $u(x^{1/r})=x,$ and thus the second integral can be bounded from above by $c|x-z|$ for some constant $c$ depending only on $s.$ Therefore, it only remains to estimate the first integral on the right hand side of \eqref{eq:gy:Hplus Lipschitz 0} in terms of $|x-z|$. For this purpose we observe that, by definition of $u$,
$$
  \left|u(x^{{1}/{r}}
  +y)-u(z^{{1}/{r}}+y)\right|
  =\left\{
  \begin{array}{ll}
  (z^{{1}/{r}}+y)^r-(x^{{1}/{r}}+y)^r
  & \text{ if }\frac{1}{2}<y<1-z^{{1}/{r}},
  \smallskip
  \\
  1-(x^{{1}/{r}}+y)^r& \text{ if }
  1-z^{{1}/{r}}<y<1-x^{{1}/{r}},
  \smallskip\\
  0 & \text{ if } 1-x^{{1}/{r}}<y.
   \end{array}\right.
$$
Therefore,
\begin{equation}
  \label{eq:gy:Hplus Lipschitz 1}
  \left|H_+(x^{{1}/{r}})-
  H_+(z^{{1}/{r}})\right|\leq
  V_1(x,z)+V_2(x,z) +c|x-z|,
  \end{equation}
where
\begin{equation}
\begin{split}
  V_1(x,z)&:=\int_{\frac{1}{2}}^{1-z^{{1}/{r}}}dy\,
  \big((z^{{1}/{r}}+y)^r-(x^{{1}/{r}}+y)^r
  \big)\frac{1}{y^{1+2s}},
  \\
  V_2(x,z)&:=\int_{1-z^{{1}/{r}}}^{1-x^{{1}/{r}}}dy\,
  \big(1-(x^{{1}/{r}}+y)^r
  \big)\frac{1}{y^{1+2s}}.
\end{split}
\end{equation}
Concerning the integrand of $V_2(x,z)$, observe that 
$1-(x^{1/r}+y)^r\leq (z^{1/r}+y)^r-(x^{1/r}+y)^r$ if $1-z^{1/r}<y.$ Hence we obtain
\begin{equation}\label{eq:gy:Hplus Lipschitz 2}
   V_1(x,z)+V_2(x,z)\leq \int_{\frac{1}{2}}^{+\infty}dy\,\big((z^{1/r}+y)^r-(x^{1/r}+y)^r\big)\frac{1}{y^{1+2s}}.
\end{equation}
By Lemma \ref{lemma:gy:technical lemma for H Lip} $(i)$ we have that
\begin{equation}
   (z^{1/r}+y)^r-(x^{1/r}+y)^r\leq 2^r(z-x)\quad\text{if }0\leq x < z\leq 2^{-r} \text{ and }
   z^{{1}/{r}}\leq y.
\end{equation}
From this,  \eqref{eq:gy:Hplus Lipschitz 1}, and 
 \eqref{eq:gy:Hplus Lipschitz 2}, it follows that
 $$
     \left|H_+(x^{{1}/{r}})-
  H_+(z^{{1}/{r}})\right|\leq
  2^r|x-z|\int_{\frac{1}{2}}^{+\infty}\frac{dy}{y^{1+2s}} +c|x-z|=c'|x-z|
$$
for some constant $c'$ depending only on $r$ and $s.$ That is, $x\mapsto H_+(x^{1/r})$ is Lipschitz in $[0,2^{-r}]$. 

To deal with $H_-$ we can proceed in a similar way as for $H_+$. We again assume \eqref{eq:gy:x smaller z} and obtain as in \eqref{eq:gy:Hplus Lipschitz 0} that (applying also a change of variables $y\mapsto -y$)
\begin{align}
  \left|H_-(x^{{1}/{r}})-
  H_-(z^{{1}/{r}})\right|
  &\leq
  \int_{\frac{1}{2}}^{+\infty}dy\,
  \frac{\big|u(x^{{1}/{r}}
  -y)-u(z^{{1}/{r}}-y)\big|}
  {y^{1+2s}}
  +
  \int_{\frac{1}{2}}^{+\infty}dy\,
  \frac{\big|u(x^{{1}/{r}}
  )-u(z^{{1}/{r}})\big|}
  {y^{1+2s}}
  \\
  &\leq \int_{\frac{1}{2}}^{+\infty}dy\,
  \frac{\big|u(x^{{1}/{r}}
  -y)-u(z^{{1}/{r}}-y)\big|}
  {y^{1+2s}}
  +
  c|x-z|.
\end{align}
This time we have (observe that now all arguments of $u$ are negative)
$$
  \left|u(x^{{1}/{r}}
  -y)-u(z^{{1}/{r}}-y)\right|
  =\left\{
  \begin{array}{ll}
  (y-x^{{1}/{r}})^r-(y-z^{{1}/{r}})^r
  & \text{ if }\frac{1}{2}<y<1+x^{{1}/{r}},
  \smallskip
  \\
1-(y-z^{{1}/{r}})^r& \text{ if }
  1+x^{{1}/{r}}<y<1+z^{{1}/{r}},
  \smallskip\\
  0 & \text{ if } 1+z^{{1}/{r}}<y.
   \end{array}\right.
$$
In the second case we have that
$1-(y-z^{1/r})^r\leq (y-x^{1/r})^r-(y-z^{1/r})^r.$ Hence we can estimate
$$
  \left|u(x^{{1}/{r}}
  -y)-u(z^{{1}/{r}}-y)\right|
  \leq
  (y-x^{{1}/{r}})^r-(y-z^{{1}/{r}})^r
  \quad\text{ for all }y>\frac{1}{2}.
$$ 
In view of Lemma \ref{lemma:gy:technical lemma for H Lip} 
$(ii)$, and using also that $0<r<1$, we now obtain 
$$
  \int_{\frac{1}{2}}^{+\infty}dy\,
  \frac{\big|u(x^{{1}/{r}}
  -y)-u(z^{{1}/{r}}-y)\big|}
  {y^{1+2s}}
  \leq
  c|x-z|\int_{\frac{1}{2}}^{+\infty}\frac{dy}{\big(y-\frac{1}{2}\big)^{1-r} y^{1+2s}}=c'|x-z|.
$$
This shows that $x\mapsto H_-(x^{1/r})$ is Lipschitz in $[0,2^{-r}]$, and concludes the proof of \eqref{eq:gy:H property}.
\end{proof}

We now prove the main result of this section. In view of its generalization in Section \ref{section:periodic counterexample i}, we formulate a lemma which is slightly more general than Proposition \ref{prop:gy:counterexample reg}.

\begin{lemma}
\label{gy:lemma:v more gener odd}
Given $0<s<1$ and $2s<r<1+2s$, let $v$ and $H$ be as in  \Cref{lemma:gy:fracLaplace of u}, and suppose that $v$
additionally satisfies
that it is odd and that
$t\mapsto H(t^{{1}/{r}})$ is Lipschitz in $[0,2^{-r}]$.
Then, there exists $f\in C^{\beta}([-2^{-r},2^{-r}])$ with $\beta=(r-2s)/r$ such that $$(-\Delta)^sv(x)=f(v(x))\quad\text{ for all ${\textstyle x\in [-\frac{1}{2},\frac{1}{2}]}$}.$$
\end{lemma}

\begin{proof}
Let us first prove the existence of $f_*\in C^{\beta}([0,2^{-r}])$ such that
\begin{equation}
\label{eq:gy:Laplace s u is fu in half interval}
  (-\Delta)^sv(x)=f_*(x^r)=f_*(v(x))\quad\text{ for all }x\in 
  {\textstyle [0,\frac{1}{2}]}.
\end{equation}  
Let $G$ and $H$ be the two functions defined in Lemma \ref{lemma:gy:fracLaplace of u}. In view of \eqref{eq:gy:fracLaplace of u}, we define 
$$
  f_*(t):=c_1 t^{\frac{r-2s}{r}}+c_2 t+c_sG(t^{{1}/{r}})
  -c_sH(t^{{1}/{r}})\quad\text{ for }t\in [0,2^{-r}].
$$
From \eqref{eq:gy:fracLaplace of u} we know that \eqref{eq:gy:Laplace s u is fu in half interval} holds for all 
$x\in {\textstyle (0,\frac{1}{2}]}$. To see that  \eqref{eq:gy:Laplace s u is fu in half interval} also holds for $x=0$, simply note that $G(0)=0$, and that $H(0)=0$ and 
$(-\Delta)^sv(0)=0$ because $v$ is odd and $v(0)=0$. 
Then, using that $2s<r$, we get 
$(-\Delta)^sv(0)=0=f_*(0)=f_*(v(0))$. 

It only remains to show that $f_*\in C^{\beta}([0,2^{-r}])$. 
Clearly $t\mapsto t^{(r-2s)/r}$ belongs to 
$C^{\beta}([0,2^{-r}]).$ Moreover, $r<1+2s$ yields $1/r>(r-2s)/r=\beta$, and hence $t\mapsto G(t^{1/r})$ belongs to $C^\beta([0,2^{-r}])$ by \eqref{eq:gy:G property}. Finally, $t\mapsto H(t^{1/r})$ is Lipschitz by assumption and, since $\beta<1$, it belongs to $C^{\beta}([0,2^{-r}])$ too. Therefore,
$f_*\in C^{\beta}([0,2^{-r}])$. 

To obtain $f$ we now take the odd extension of $f_*,$ i.e., 
\begin{equation}
\text{$f(x):=-f_*(-x)$\quad if $x\in [-2^{-r},0)$\qquad and\qquad  $f(x):=f_*(x)$\quad if $x\in [0,2^{-r}].$}
\end{equation} 
Since $f_*\in C^{\beta}([0,2^{-r}])$ and $f_*(0)=0$, we then get that
$f\in C^{\beta}([-2^{-r},2^{-r}])$.

It is easy check that $(-\Delta)^s v$ is an odd function because $v$ is an odd function. Using \eqref{eq:gy:Laplace s u is fu in half interval}, this yields for every $x\in {\textstyle [0,\frac{1}{2}]}$ that 
$$
  (-\Delta)^sv(-x)=-(-\Delta)^sv(x)=-f_*(v(x))=f(-v(x))=f(v(-x)).
$$
Hence $(-\Delta)^sv(x)=f(v(x))$ for all ${\textstyle x\in [-\frac{1}{2},\frac{1}{2}]}$ and this concludes the proof of the lemma.
\end{proof}

\smallskip

We now give the proof of  \Cref{prop:gy:counterexample reg}, which follows from \Cref{lemma:gy:fracLaplace of u,gy:lemma:v more gener odd}.

\begin{proof}[Proof of Proposition \ref{prop:gy:counterexample reg}.]
Let $u$ be the function defined in \eqref{eq:gy:u def} with 
$$
  r=\frac{2s}{1-\beta}, \quad\text{which is equivalent to }\beta=\frac{r-2s}{r}.
$$
Clearly $u\in C^r(\R)$ (note that $u$ is indeed Lipschitz away from the origin and $2s/(1-\beta)<1$ by assumption) but 
$u\notin C^{r+\epsilon}([-\rho,\rho])$ for any $\epsilon,\rho>0.$ The fact that $(-\Delta)^su=f(u)$ in ${\textstyle [-\frac{1}{2},\frac{1}{2}]}$ for some function $f\in C^{\beta}(\R)$ follows from \Cref{lemma:gy:fracLaplace of u} (note that $r<1$ because $2s<1-\beta$ by assumption) and \Cref{gy:lemma:v more gener odd} applied to $u$.
\end{proof}

\smallskip

We now give an example that shows the optimality of the exponent $2s/(1-\beta)$ in \Cref{fthm:Holder reg periodic}~$(i)$ (in the nonperiodic setting) for the case $2s=1-\beta.$ This example is much simpler than the one presented above for the case $2s<1-\beta$, and also not very surprising as in this case $2s/(1-\beta)=1=2s+\beta$ is indeed the expected gain of regularity for solutions to the linear equation 
$(-\Delta)^su=g$ with $g$ of class $C^{\beta}$. 
In this case we define $u$ by
\begin{equation}
 \label{gy:eq:u is modulus x and even}
  u(x):=\left\{\begin{array}{rl}
        |x|& \text{ if }|x|\leq 1,
        \\
        1& \text{ if }|x|> 1.
        \end{array}\right.
\end{equation}
Clearly $u\notin C^{1+\epsilon}(\R)=C^{2s/(1-\beta)+\epsilon}(\R)$ for any $\epsilon>0$, but $u$ is a bounded solution to the following semilinear equation stated in the next lemma.

\begin{lemma}
\label{lemma:regularity example no surprise} Let $0<s<1$ and $\beta>0$ such that
$2s=1-\beta$, and let $u$ be defined by \eqref{gy:eq:u is modulus x and even}. Then, 
$$
  (-\Delta)^s u(x)=f(u(x))\quad\text{ for all }{\textstyle x\in [-\frac{1}{2},\frac{1}{2}]},
$$
for some measurable function $f:\R\to\R$ such that $f\in C^{\beta}([0,{\textstyle \frac{1}{2}}]).$ 
\end{lemma}

\begin{proof}
For this $u$ we can compute $(-\Delta)^su$ explicitly.
Given $x\in [0,{\textstyle \frac{1}{2}}]$, note that $u(x)=x$ and 
$$
  u(x+y)=\left\{\begin{array}{rl}
         1 & \text{ if }y\in(-\infty,-1-x)=:Q_1,
         \\
         -(x+y) & \text{ if }y\in (-1-x,-x)=:Q_2,
         \\
         x+y & \text{ if }y\in(-x,1-x)=:Q_3,
         \\
         1 & \text{ if }y\in(1-x,+\infty)=:Q_4\,.
         \end{array}\right.
$$
Therefore, recalling that $1-2s=\beta>0$ by assumption, we have
\begin{align}
  \int_{Q_1}dy\,\frac{u(x+y)-u(x)}{|y|^{1+2s}}
  &=
  \frac{1-x}{2s(1+x)^{2s}},
  \\
  \int_{Q_2}dy\,\frac{u(x+y)-u(x)}{|y|^{1+2s}}
  &=
  \frac{1}{s}\big(x(1+x)^{-2s}-x^{1-2s}\big)
  -\frac{1}{1-2s}\big(x^{1-2s}-(1+x)^{1-2s}\big),
  \\
  \int_{Q_3}dy\,\frac{u(x+y)-u(x)}{|y|^{1+2s}}
  &=
  \frac{1}{1-2s}\big((1-x)^{1-2s}-x^{1-2s}\big),
  \\
  \int_{Q_4}dy\,\frac{u(x+y)-u(x)}{|y|^{1+2s}}
  &=
  \frac{1}{2s}(1-x)^{1-2s}.
\end{align}
Collecting the terms containing $x^{1-2s}$ and powers of $(1\pm x)$ or $x(1+x)^{-2s}$ gives that 
\begin{equation}\label{gy:example.2s=1-beta.eq1}
  (-\Delta)^su(x)=c_s\frac{x^{1-2s}}{s(1-2s)}+G(x)=:f(x)\quad\text{for all }x\in [0,{\textstyle \frac{1}{2}}],
\end{equation}
for some $G\in C^{\infty}([0,{\textstyle \frac{1}{2}}])$.
In particular, $f\in C^{\beta}([0,{\textstyle \frac{1}{2}}])$ but 
$f\not\in C^{\beta+\epsilon}([0,\rho])$ for any $\epsilon,\rho>0$. Now, for $x\in [0,1]$ we have $x=u(x)$ and, thus,
$(-\Delta)^s u(x)=f(u(x)).$ It is easy to see that $(-\Delta)^su$ is an even function because $u$ is an even function. We therefore obtain that
$$
  (-\Delta)^s u(x)=f(u(x))\quad\text{ for all }{\textstyle x\in [-\frac{1}{2},\frac{1}{2}]},
$$
which concludes the proof of the lemma.
\end{proof}

\section{The periodic example for $2s\leq 1-\beta$}
\label{section:periodic counterexample i}

In this section we adapt the examples of \Cref{section:nonperiodic counterexample i} (which dealt with the nonperiodic setting) to give the proof of \Cref{gy:thm:intro:examples reg} $(i)$. 

\begin{proof}[Proof of \Cref{gy:thm:intro:examples reg} $(i)$]
Since $(-\Delta)^s (u(\lambda \cdot))(x)
=\lambda^{2s}(-\Delta)^s u(\lambda x)$ for all $\lambda>0$,  the proof follows from the following proposition.
\end{proof}

\begin{proposition}\label{gy:example_8periodic:lemma}
Let $0<\beta<1$ and $0<2s\leq 1-\beta.$ There exist two functions $u\in C^{\frac{2s}{1-\beta}}(\R)$ and $f\in C^{\beta}(\R)$ such that $u$ is periodic, $(-\Delta)^s u=f(u)$ in $\R,$ and
$$
  u\notin C^{\frac{2s}{1-\beta}+\epsilon}([-\rho,\rho])\quad\text{for any $\epsilon,\rho>0$}.
$$
\end{proposition}

To prove this proposition we will adapt the function given in  \eqref{eq:gy:u def} to the periodic setting. We will only deal with the case $2s<1-\beta$ in detail; the case $2s=1-\beta$ can be dealt with similarly by taking a periodized version of the function given in  \eqref{gy:eq:u is modulus x and even}. 
We construct the function $u$ in \Cref{gy:example_8periodic:lemma} to be $8$-periodic, and chosen to satisfy the following properties:
\begin{itemize}
 \item[$(a)$] Set $r:=2s/(1-\beta)$ and
 $$ 
   u(x):=\left\{\begin{array}{rl}
        -|x|^r & \text{ if }-1<x<0,
        \\
        x^r & \text{ if }0\leq x<1.
        \end{array}\right.
 $$ 
 
 \item[$(b)$] $u$ is smooth in $(-4,0)\cup(0,4)$.
 
 \item[$(c)$] $u$ is increasing in $(0,2).$

 \item[$(d)$] $u$ is odd with respect to $x=0$ and even with respect to $x=2$.

 \item[$(e)$]
 The graph of $u$ is flatter in $(1,2)$ than in $(0,1)$. That is,
 $$
   \sup\{u'(x):\, x\in(1,2)\}= u'(1)=\inf\{u'(x):\, x\in(0,1)\}.
 $$
 \item[$(f)$] $u$ is a quadratic function near $x=2$. More precisely, there exists $\delta>0$ such that 
 $$
   u(x)=-(x-2)^2+u(2)\qquad\text{for all }x\in (2-\delta,2+\delta).
 $$
\end{itemize}
An immediate consequence of $(d)$ is that 
\begin{equation}\label{gy:eq.x+4to-x}
  u(x+4)=-u(x)\quad\text{for all }x\in\R.
\end{equation}

The main ingredient of the proof of \Cref{gy:example_8periodic:lemma} is again  \Cref{lemma:gy:fracLaplace of u}, but we have to show  \eqref{eq:gy:H property} for the new periodic function $u.$ This is the content of the next lemma.

\begin{lemma}
\label{lemma:gy:H property u periodic}
Let $u$ be defined by $(a)$-$(f)$ and $H$ be defined by
$$
  H(x):=\int_{\R\setminus(-\frac{1}{2},\frac{1}{2})}
  dy\,\frac{u(x+y)-u(x)}{|y|^{1+2s}}.
$$
Then, $t\mapsto H(t^{{1}/{r}})$ is  Lipschitz in $[0, 2^{-r}]$.
\end{lemma}

\begin{proof} We decompose $H$ as in \eqref{eq:gy:decomp of H} into $H_+$ and $H_-$ by splitting the range of integration into $(\frac{1}{2},+\infty)$ and $(-\infty,-\frac{1}{2})$.
Let us first show that $t\mapsto H_+(t^{1/r})$ is Lipschitz in $[0,2^{-r}].$ Assume that
\begin{equation}\label{eq:gy:x smaller z periodic case}
   0\leq x<z\leq 2^{-r},\quad\text{and thus }\quad
  0\leq x^{{1}/{r}}<z^{{1}/{r}}\leq\frac{1}{2}.
\end{equation}
We obviously have that 
\begin{align}
    \left|H_+(x^{{1}/{r}})-
  H_+(z^{{1}/{r}})\right|
  \leq &
  \int_{\frac{1}{2}}^{+\infty}dy\,
  \frac{\big|u(x^{{1}/{r}}
  +y)-u(z^{{1}/{r}}+y)\big|}
  {y^{1+2s}}
  +
  \int_{\frac{1}{2}}^{+\infty}dy\,
  \frac{\big|u(x^{{1}/{r}}
  )-u(z^{{1}/{r}})\big|}
  {y^{1+2s}}
  \\
  \leq & 
  \int_{\frac{1}{2}}^{+\infty}dy\,
  \frac{\big|u(x^{{1}/{r}}
  +y)-u(z^{{1}/{r}}+y)\big|}
  {y^{1+2s}}
  + C(s)|x-z|.
\end{align}
Therefore, it is sufficient to prove that
\begin{equation}\label{eq:gy:main estim for H periodic case}
  \int_{\frac{1}{2}}^{+\infty}dy\,
  \frac{\big|u(x^{{1}/{r}}
  +y)-u(z^{{1}/{r}}+y)\big|}
  {y^{1+2s}}
  \leq
  C(r,s)|x-z|.
\end{equation}
Due to the piecewise definition of $u,$ we split the interval $(\frac{1}{2},+\infty)$ in the following way: For $k=0,1,2,\ldots$ define
\begin{align}
  A_1^k&:=\left({\textstyle\frac{1}{2}},4-z^{{1}/{r}}-x^{{1}/{r}}\right)
  +8k,
  \\
  A_2^k&:=\left(4-z^{{1}/{r}}-x^{{1}/{r}},4-z^{{1}/{r}}  
  \right)
  +8k,
  \\
  A_3^k&:=\left(4-z^{{1}/{r}},4-x^{{1}/{r}}  
  \right)
  +8k,
  \\
  A_4^k&:=\left(4-x^{{1}/{r}},4+{\textstyle\frac{1}{2}}  
  \right)
  +8k,
  \\
  B_i^k&:= A_i^k+4\quad\text{for }i=1,2,3,4. 
\end{align}
With these definitions we see that $({\textstyle\frac{1}{2}},+\infty)$ coincides with the disjoint union
$$
  \left({\textstyle\frac{1}{2}},+\infty\right)=\bigcup_{k=0}^{+\infty}
  \bigcup_{i=1}^4 \left(A_i^k\cup B_i^k\right)
  \quad\text{ up to a set of measure zero.}
$$
To prove the inequality \eqref{eq:gy:main estim for H periodic case} we will now show that
\begin{equation}
 \label{eq:gy:Aik estimates Lip}
   \sum_{k=0}^{+\infty}
     \int_{A_i^k}dy\,
  \frac{\big|u(x^{{1}/{r}}
  +y)-u(z^{{1}/{r}}+y)\big|}
  {y^{1+2s}}
  \leq
  C(r,s)|x-z|\qquad\text{for }i=1,2,3,4.
\end{equation}
The corresponding result for $B_i^k$, i.e. \eqref{eq:gy:Aik estimates Lip} with $A_i^k$ replaced by $B_i^k,$ follows then immediately from \eqref{eq:gy:Aik estimates Lip} using that $u(t+4)=-u(t)$ for all $t\in\R$. 

\medskip
\noindent\underline{{\em Estimate for $A_1^k$}\,:} 
\smallskip

Let us firs estimate $|u(x^{1/r}+y)-u(z^{1/r}+y)|.$ We claim that
\begin{equation}
 \label{eq:gy:A1k estimate}
  \big|u(x^{{1}/{r}}
  +y)-u(z^{{1}/{r}}+y)\big|
  \leq
  \big|u(x^{{1}/{r}}
  )-u(z^{{1}/{r}})\big|=|x-z|
  \quad\text{if }y\in A_1^k.
\end{equation}
Once \eqref{eq:gy:A1k estimate} is proven it follows that, for $i=1$, the left hand side of \eqref{eq:gy:Aik estimates Lip} is smaller than
$$
  |x-z|\sum_{k=0}^{+\infty}\int_{A_1^k}\frac{dy}{y^{1+2s}}
  \leq
  |x-z|\int_{\frac{1}{2}}^{+\infty}\frac{dy}{y^{1+2s}}=C(s)|x-z|.
$$

We now prove \eqref{eq:gy:A1k estimate}.
By the periodicity of $u$ we can assume without loss of generality that $k=0.$ We now use the properties $(a)$, $(d)$, and $(e)$ of $u$ which give that the graph of $u$ is flatter in 
$(1,3)$ than in $(0,1)$ or $(3,4)$, and that the graph of $u$ in $(0,4)$ is steeper the closer we are to $0$ or $4$ (recall that $0<r<1$). Therefore, geometrically \eqref{eq:gy:A1k estimate} is obvious (observe also that $x^{1/r}$ is not further from the origin than $z^{1/r}+y$ from 4), but nevertheless we give a proof of it. 

Consider first the simpler case when $x^{1/r}+y,z^{1/r}+y\in (\frac{1}{2},3).$ Then, using that
\begin{equation}\label{gy:eq.aux.A02}
  \sup_{t\in[\frac{1}{2},3)}|u'(t)|=u'({\textstyle\frac{1}{2}})
  =\inf_{t\in (0,\frac{1}{2}]}u'(t)
\end{equation}
and that $x^{{1}/{r}},z^{{1}/{r}}\in[0,\frac{1}{2}]$, one easily obtains \eqref{eq:gy:A1k estimate}.

In the case that 
$z^{1/r}+y\in [3,4-x^{1/r})$ we have $x^{1/r}+y\in (2,4-z^{1/r})$ because $|x^{{1}/{r}}-z^{{1}/{r}}|\leq\frac{1}{2}$. This leads to 
\begin{equation}\label{gy:eq.aux.A03}
4-z^{{1}/{r}}-y\in(x^{1/r},1]\quad\text{and}\quad
4-x^{{1}/{r}}-y\in(z^{1/r},2).
\end{equation}
Then, using \eqref{gy:eq.aux.A02} in the case that 
$1\leq4-x^{{1}/{r}}-y<2$ (which yields 
$\frac{1}{2}\leq4-z^{{1}/{r}}-y\leq1$), or using $(a)$ in the case that 
$z^{1/r}<4-x^{{1}/{r}}-y<1$ (which gives $x^{1/r}<4-z^{{1}/{r}}-y<1$), we get
\begin{equation}\label{gy:eq.aux.A01}
\big|u(z^{{1}/{r}})-u(x^{{1}/{r}})\big|\geq
\big|u(4-x^{{1}/{r}}-y)-u(4-z^{{1}/{r}}-y)\big|
=\big|u(x^{{1}/{r}}+y)-u(z^{{1}/{r}}+y)\big|;
\end{equation}
the last equality above follows from that fact that $u$ is even with respect to $2$. Therefore, \eqref{eq:gy:A1k estimate} also holds in the later case $z^{1/r}+y\in [3,4-x^{1/r})$.

\medskip
\noindent\underline{{\em Estimate for $A_2^k$}\,:} 
\smallskip

First, use the inclusion $A_2^k\subset(3,4-z^{1/r})+8k$ to get
$$
     \int_{A_2^k}dy\,
  \frac{\big|u(x^{{1}/{r}}
  +y)-u(z^{{1}/{r}}+y)\big|}
  {y^{1+2s}}
  \leq
  \int_{3+8k}^{4-z^{{1}/{r}}+8k}dy\,
  \frac{\big|u(x^{{1}/{r}}
  +y)-u(z^{{1}/{r}}+y)\big|}
  {y^{1+2s}}.
$$
We will apply the change of  variables 
$y=-y'+8k+4$ with $y'\in (z^{1/r},1)$. Observe that then 
$x^{1/r}-y',z^{1/r}-y'\in(-1,0)$, and thus
$$
  u(x^{{1}/{r}}+y)= u(x^{{1}/{r}}-y'+8k+4)
  =u(x^{{1}/{r}}-y'+4)
  =-u(x^{{1}/{r}}-y')
  =(y'-x^{{1}/{r}})^r.
$$
By the same calculation we also obtain that
$
  u(z^{{1}/{r}}+y)
  =(y'-z^{{1}/{r}})^r.
$
From this, and using also 
\Cref{lemma:gy:technical lemma for H Lip}~$(ii)$, we get
\begin{equation}
\begin{split}
   \int_{A_2^k}dy\,
  \frac{\big|u(x^{{1}/{r}}
  +y)-u(z^{{1}/{r}}+y)\big|}
  {y^{1+2s}}
  &\leq 
   \int_{z^{{1}/{r}}}^{1}dy'\,
  \frac{(y'-x^{{1}/{r}})^r-(y'-z^{{1}/{r}})^r}
  {(8k+4-y')^{1+2s}}
  \\
  &\leq 
  \frac{2^{r-1}(z-x)}{(8k+3)^{1+2s}}\int_{z^{{1}/{r}}}^1
  \frac{dy'}{(y'-z^{{1}/{r}})^{1-r}}
  \\
  &\leq
  \frac{2^{r-1}(z-x)}{(8k+3)^{1+2s}}
  \int_0^1 \frac{dt}{t^{1-r}}=\frac{C(r)|x-z|}{(8k+3)^{1+2s}}.
\end{split}
\end{equation}
Therefore, 
$$
     \sum_{k=0}^{+\infty}
     \int_{A_2^k}dy\,
  \frac{\big|u(x^{{1}/{r}}
  +y)-u(z^{{1}/{r}}+y)\big|}
  {y^{1+2s}}
  \leq
  C(r)|x-z|\sum_{k=0}^{+\infty}\frac{1}{(8k+3)^{1+2s}}
  =C(r,s)|x-z|.
$$

\medskip
\noindent\underline{{\em Estimate for $A_3^k$}\,:} 
\smallskip

For $y\in A_3^k=(-z^{1/r}+8k+4,-x^{1/r}+8k+4)$ we again do the change of variables $y=-y'+8k+4$, now with 
$y'\in (x^{1/r},z^{1/r}).$ This time we have
$x^{1/r}-y'\in (-1,0)$ but $z^{1/r}-y'\in(0,1),$ which yields
\begin{align}
  u(x^{{1}/{r}}+y)&=u(x^{{1}/{r}}-y'+8k+4)=-u(x^{{1}/{r}}-y')=(y'-x^{{1}/{r}})^r,
  \\
  u(z^{{1}/{r}}+y)&=u(z^{{1}/{r}}-y'+8k+4)=-u(z^{{1}/{r}}-y')=-(z^{{1}/{r}}-y')^r.
\end{align}
Hence we obtain
\begin{align}
  \int_{A_3^k}dy\,&
  \frac{\big|u(x^{{1}/{r}}
  +y)-u(z^{{1}/{r}}+y)\big|}
  {y^{1+2s}}
  = 
  \int_{x^{{1}/{r}}}^{z^{{1}/{r}}}dy'\,\frac{(y'-x^{{1}/{r}})^r+ 
  (z^{{1}/{r}}-y')^r}{(8k+4-y')^{1+2s}} 
  \\
  &\quad \leq
  \frac{1}{(8k+3)^{1+2s}} 
  \int_{x^{{1}/{r}}}^{z^{{1}/{r}}}dy'\,
  \big((y'-x^{{1}/{r}})^r+ 
  (z^{{1}/{r}}-y')^r\big)
  =
  \frac{2(z^{{1}/{r}}-x^{{1}/{r}})^{r+1}}{(r+1)(8k+3)^{1+2s}}
  \\
  &\quad 
  =
  \frac{2}{(r+1)(8k+3)^{1+2s}}\,
  (z^{{1}/{r}}-x^{{1}/{r}})^{r}\,
  \frac{z^{{1}/{r}}-x^{{1}/{r}}}{z-x}\, (z-x).
\end{align}

Proceeding as in the proof of \Cref{lemma:gy:technical lemma for H Lip}~$(ii)$ and using that $z^{1/r}-x^{1/r}\leq\frac{1}{2}$ we get
$$
   \int_{A_3^k}dy\,
  \frac{\big|u(x^{{1}/{r}}
  +y)-u(z^{{1}/{r}}+y)\big|}
  {y^{1+2s}}
  \leq \frac{2}{(r+1)(8k+3)^{1+2s}}\, 2^{-r}\,\frac{2^{r-1}}{r}\,(z-x).
$$
Thus, we can conclude exactly as in the last step of the estimate for $A_2^k$, summing up in $k$.

\medskip
\noindent\underline{{\em Estimate for $A_4^k$}\,:} 
\smallskip

Now $y\in (4-x^{1/r}+8k,4+\frac{1}{2}+8k).$  This time we do the change of variables $y=y'+8k+4$ with $y'\in (-x^{1/r},\frac{1}{2}).$  Hence $x^{1/r}+y',z^{1/r}+y'\in (0,1)$ and we obtain
\begin{align}
  u(x^{{1}/{r}}+y)&=u(x^{{1}/{r}}+y'+8k+4)=-u(x^{{1}/{r}}+y')
  =-(x^{{1}/{r}}+y')^r,
  \\
  u(z^{{1}/{r}}+y)&=u(z^{{1}/{r}}+y'+8k+4)=-u(z^{{1}/{r}}+y')=-(z^{{1}/{r}}+y')^r.
\end{align}
This leads to
\begin{align}
  \int_{A_4^k}dy\,
  \frac{\big|u(x^{{1}/{r}}
  +y)-u(z^{{1}/{r}}+y)\big|}
  {y^{1+2s}}
  \leq 
  \frac{1}{(8k+3)^{1+2s}}\int_{-x^{{1}/{r}}}^{\frac{1}{2}}dy'\,
  \left((z^{{1}/{r}}+y')^r
  -(x^{{1}/{r}}+y')^r\right).
\end{align}
We proceed similarly as in the proof of  \Cref{lemma:gy:technical lemma for H Lip} $(ii)$ and obtain that the integrand can be estimated by 
$$
  \frac{(z^{{1}/{r}}+y')^r
  -(x^{{1}/{r}}+y')^r}{z^{{1}/{r}}-x^{{1}/{r}}}
  \, \frac{z^{{1}/{r}}-x^{{1}/{r}}}{z-x}\,(z-x)
  \leq
  \varphi'(t_0)\frac{2^{r-1}}{r}(z-x),
$$
where $\varphi(t):=(t+y')^r,$ $t_0\in(x^{1/r},z^{1/r})$, and 
$\varphi'(y_0)\leq r (x^{1/r}+y')^{r-1}.$ We can therefore estimate
\begin{align}
  \int_{A_4^k}dy\,
  \frac{\big|u(x^{{1}/{r}}
  +y)-u(z^{{1}/{r}}+y)\big|}
  {y^{1+2s}}
  &\leq 
  \frac{2^{r-1}(z-x)}{(8k+3)^{1+2s}}\int_{-x^{{1}/{r}}}
  ^{\frac{1}{2}}\frac{dy'}{(y'+x^{{1}/{r}})^{1-r}}
  \\
  &\leq  
  \frac{2^{r-1}(z-x)}{(8k+3)^{1+2s}}\int_0^1 \frac{dt}{t^{1-r}}
  =\frac{C(r)|x-z|}{(8k+3)^{1+2s}}.
\end{align}
We can now conclude as in the estimate for $A_2^k$ or $A_3^k$ and this gives the proof of \eqref{eq:gy:Aik estimates Lip} for $A_4^k$.
\smallskip

The four estimates in \eqref{eq:gy:Aik estimates Lip} for each $A_i^k$, $i=1,2,3,4,$ yield \eqref{eq:gy:main estim for H periodic case} and show that $t\mapsto H_+(t^{1/r})$ is Lipschitz in $[0, 2^{-r}]$.
The corresponding result for $H_-$ ---recall its definition in \eqref{eq:gy:decomp of H}--- can be easily deduced from  what we have shown for $H_+$. As for $H_+$, it is sufficient to prove \eqref{eq:gy:main estim for H periodic case} but with the integration over $(\frac{1}{2},+\infty)$ replaced by $(-\infty,-\frac{1}{2})$ and $y$ by $|y|$. Now, observe that (doing the change of variables $y=-y'+16$ in the third line)  
\begin{align}
   \int_{-\infty}^{-\frac{1}{2}}dy\,&
  \frac{\big|u(x^{{1}/{r}}
  +y)-u(z^{{1}/{r}}+y)\big|}
  {|y|^{1+2s}}
  =
  \sum_{k=0}^{+\infty}\int_{\frac{1}{2}+8k}^{\frac{1}{2}+8k+8}
  dy\,
  \frac{\big|u(x^{{1}/{r}}
  -y)-u(z^{{1}/{r}}-y)\big|}
  {y^{1+2s}}
  \\
  &\qquad\leq \sum_{k=0}^{+\infty}\frac{1}{\left(\frac{1}{2}+8k\right)^{1+2s}}
  \int_{\frac{1}{2}}^{\frac{1}{2}+8}dy\,|u(x^{{1}/{r}}
  -y)-u(z^{{1}/{r}}-y)\big|
  \\
  &\qquad=
  \sum_{k=0}^{+\infty}\frac{1}{\left(\frac{1}{2}+8k\right)^{1+2s}}
  \int_{8-\frac{1}{2}}^{16-\frac{1}{2}}dy'\,|u(x^{{1}/{r}}
  +y')-u(z^{{1}/{r}}+y')\big|
  \\
  &\qquad\leq 
  \sum_{k=0}^{+\infty}\frac{\left(16-\frac{1}{2}+8k\right)^{1+2s}}{\left(\frac{1}{2}+8k\right)^{1+2s}}
  \int_{8-\frac{1}{2}}^{16-\frac{1}{2}}dy'\,\frac{|u(x^{{1}/{r}}
  +y')-u(z^{{1}/{r}}+y')\big|}{|y'+8k|^{1+2s}}
  \\
  &\qquad\leq 
  \frac{\left(16-\frac{1}{2}\right)^{1+2s}}{\left(\frac{1}{2}\right)^{1+2s}}
  \sum_{k=0}^{+\infty}
  \int_{8-\frac{1}{2}}^{16-\frac{1}{2}}dy'\,\frac{|u(x^{{1}/{r}}
  +y')-u(z^{{1}/{r}}+y')\big|}{|y'+8k|^{1+2s}}
  \\
  &\qquad\leq  32^3 \sum_{k=0}^{+\infty}
  \int_{8-\frac{1}{2}}^{16-\frac{1}{2}}dy'\,\frac{|u(x^{{1}/{r}}
  +y'+8k)-u(z^{{1}/{r}}+y'+8k)\big|}{|y'+8k|^{1+2s}}
  \\
  &\qquad=
  32^3
  \int_{8-\frac{1}{2}}^{+\infty}dy\,
  \frac{\big|u(x^{{1}/{r}}
  +y)-u(z^{{1}/{r}}+y)\big|}
  {y^{1+2s}}
  \leq C(r,s)|x-z|,
\end{align}
where we have used \eqref{eq:gy:main estim for H periodic case} in the last inequality.
\end{proof}

We now give the proof of the main result of this section, namely \Cref{gy:example_8periodic:lemma}.
As already mentioned, we only give the proof for the case $2s<1-\beta$; see the comment below \Cref{gy:example_8periodic:lemma} for the case $2s=1-\beta$. 

\begin{proof}[Proof of \Cref{gy:example_8periodic:lemma}]
Since $0<\beta<1$, $0<2s< 1-\beta$, and $r=2s/(1-\beta)$, we have $0<2s<r<1$.
In view of  \Cref{lemma:gy:fracLaplace of u} and of \Cref{gy:lemma:v more gener odd}, whose hypothesis on the Lipschitz property of $t\mapsto H(t^{1/r})$ is satisfied by   \Cref{lemma:gy:H property u periodic}, we know that
\begin{equation}\label{gy:implicit.def.f.periodic}
  \text{there exists $f\in C^{\beta}([0,2^{-r}])$ such that }(-\Delta)^su(x)=f(u(x))\text{ for all }x\in {\textstyle [0,\frac{1}{2}]}.
\end{equation} 
Indeed, $f$ can be extended to 
$[0,\max u]=[0,u(2)]$ using
$(-\Delta)^su(x)=f(u(x))$ as an implicit equation for $f$. To see that such an $f$ exists and is unique, recall that $u$ is smooth and increasing in $(0,2)$ and, thus, the inverse $u^{-1}$ is smooth in $(0,u(2))$ and continuous in $[0,u(2)]$. This leads to define
\begin{equation}\label{gy:implicit.def.f.periodic.truedef}
  f(t):=\left((-\Delta)^s u\right)(u^{-1}(t))\quad\text{for all $t\in[0,u(2)]$.}
\end{equation} 
With this definition, we see that $f$ is continuous in $[0,u(2)]$ and smooth in $(0,u(2)),$ as $(-\Delta)^su$ is smooth in $(0,4)$. In addition, $(-\Delta)^su(x)=f(u(x))$ for all $x\in[0,2]$ by \eqref{gy:implicit.def.f.periodic.truedef}. Moreover, by \eqref{gy:implicit.def.f.periodic} and the uniqueness of $f$ (also in $[0,2^{-r}]$), we know that  
$f\in C^{\beta}([0,2^{-r}])$.

We will now show that $f\in C^{\beta}([0,u(2)])$ and, for this purpose, we will use property $(f)$ of $u$. Let us abbreviate $g(x):=(-\Delta)^su(x).$ As $u$ and $g$ are smooth in $(2-\delta,2)$ for some $\delta>0$, and $f$ is smooth in $(u(2-\delta),u(2))$ we see that
\begin{equation}\label{eq:gy:f prime}
  f'(u(x))u'(x)=g'(x)\quad\text{for all }x\in (2-\delta,2).
\end{equation}
As $u$ is even with respect to $x=2$ we obtain that $u'$ is an odd function with respect to $x=2$ and  $g'(2)=((-\Delta)^s u')(2)=0$ ---here we have also used that $u,$ and hence $g$ too, are smooth in a neighborhood of $x=2$. Using the smoothness of $g$ and a Taylor expansion we get
$$
  g'(x)=g'(2)+O(|x-2|)=O(|x-2|)\quad\text{for }x\to 2.
$$
Therefore, from \eqref{eq:gy:f prime} and property $(f)$ of $u$ we see that, for some constant $C>0$,  
$$
  |2f'(u(x))(x-2)|\leq C|x-2|\quad\text{for all }x\in (2-\delta,2).
$$ 
Hence $f'$ is bounded in $(u(2-\delta),u(2)).$ As $f$ is also continuous up to $u(2)$, we get that $f$ is Lipschitz in 
$[u(2-\delta),u(2)]$. Combining this with the fact that $f\in C^{\beta}([0,2^{-r}])$ and that $f$ is smooth in $(0,u(2))$, we conclude that $f\in C^{\beta}([0,u(2)])$.

Arguing now as in the proof of \Cref{gy:lemma:v more gener odd}, if we take the odd extension of $f$ we have that
\begin{equation}\label{eq:gy:x in minus two two}
  (-\Delta)^su(x)=f(u(x))\quad\text{ for all }x\in [-2,2],\quad\text{ where }f\in C^{\beta}([-u(2),u(2)]).
\end{equation}
It remains to show that 
\begin{equation}\label{eq:gy:extending}
(-\Delta)^su(x)=f(u(x))
\end{equation}
for all $x\in [-4,-2)\cup(2,4]$. We will only give the details for the case $x\in(2,4]$, the case $x\in[-4,-2)$ is analogous. Let $x\in (2,4]$ and $x'\in[0,2)$ be the reflection of $x$ with respect to $2,$ i.e., $x-2=2-x'$. Then, clearly, $f(u(x))=f(u(x'))$ because $u$ is even with respect to $2.$ Since $(-\Delta)^su$ is also even with respect to $2$ we obtain from \eqref{eq:gy:x in minus two two} that \eqref{eq:gy:extending} holds for all 
$x\in (2,4]$. Analogously, since $u$ is also even with respect to $-2$, we also get \eqref{eq:gy:extending} for all $x\in[-4,-2)$. 

Up to now we have seen that \eqref{eq:gy:extending} holds for all 
$x\in [-4,4]$ for some $f\in C^{\beta}([\min u,\max u])$. Now, using the periodicity of $u$, we conclude that
\eqref{eq:gy:extending} holds
for all $x\in \R$
for some $f\in C^{\beta}(\R),$ where $f$ has been extended in an arbitrary $C^{\beta}$ way from $[\min u,\max u]$ to the whole real line.
\end{proof}

\section{The case $2s> 1-\beta$}
\label{section:nonperiodic counterexample ii}

The goal of this section is to prove  \Cref{gy:thm:intro:examples reg} $(ii)$. As for the case $2s\leq 1-\beta$, we will first show the optimality of a local regularity result in a nonperiodic setting, that is, we will give the counterpart of \Cref{prop:gy:counterexample reg} and \Cref{lemma:regularity example no surprise} for 
$2s> 1-\beta$. This is the content of the following proposition. 

\begin{proposition}
\label{prop:gy:counterexample reg ii}
Let $0<s<1$ and $0<\beta<1$ such that $2s> 1-\beta$.
Then, there exists
$$
  u\in C^{2s+\beta}(\R)\quad\text{ such that }\quad u\notin C^{2s+\beta+\epsilon}([-\rho,\rho])
$$
for any $\epsilon,\rho>0$ and there exists $f\in C^{\beta}(\R)$ such that 
$$
  (-\Delta)^su(x)=f(u(x))\quad\text{for all }
  x\in {\textstyle[-\frac{1}{2},\frac{1}{2}]}.
$$
\end{proposition}

The function that will satisfy the claims of this proposition is
\begin{equation}\label{gy:eq:u for ii nonper}
  u(x):=\left\{
        \begin{array}{rl}
        \sgn(x)\,|x|^{2s+\beta}+x & \text{if }|x|\leq 1,
        \\
		2\, \sgn(x) & \text{if }|x|>1;
		\end{array}\right.
\end{equation}
see \Cref{section:appendix opt reg local} for a motivation why we have added $x$ to the power function given by \eqref{eq:gy:u def}, which was used in the case $2s< 1-\beta.$ Let us briefly say here that we need to make the inverse $u^{-1}$ more regular, since the inverse of $x^{2s+\beta}$ is not Lipschitz if $2s+\beta>1.$

\begin{proof}[Proof of \Cref{prop:gy:counterexample reg ii}]
We write the function $u$ as $u=v+\varphi,$ where
\begin{equation}\label{gy:eq:split u in v plus varphi}
  v(x):=\left\{
        \begin{array}{rl}
        x^{2s+\beta} & \text{if }0\leq x\leq1,
        \\
		-|x|^{2s+\beta} & \text{if }-1\leq x< 0,
		\\
		\sgn(x) & \text{if }|x|> 1,
		\end{array}\right.
		\qquad
	\varphi(x):=\left\{
        \begin{array}{rl}
        x & \text{if }|x|\leq1,
		\\
		\sgn(x) & \text{if }|x|> 1.
		\end{array}\right.
\end{equation}
Now, assume that $x\in [0,\frac{1}{2}]$ and note that if 
$|x+y|\leq1$ then $\varphi(x)-\varphi(x+y)=-y.$ We therefore obtain  
(the first integral is understood in the principal value sense)
\begin{align}
  (-\Delta)^s\varphi(x)&=-c_s\int_{-\frac{1}{2}}^{\frac{1}{2}}dy\, 
  \frac{y}{|y|^{1+2s}}
  +c_s\int_{\R\setminus(-\frac{1}{2},\frac{1}{2})}dy\, 
  \frac{\varphi(x)-\varphi(x+y)}{|y|^{1+2s}}
  \\
  &=c_s\int_{\R\setminus(-\frac{1}{2},\frac{1}{2})}dy\, 
  \frac{\varphi(x)-\varphi(x+y)}{|y|^{1+2s}}.
\end{align}
We then apply \Cref{lemma:gy:fracLaplace of u} with $r=2s+\beta$
to calculate $(-\Delta)^sv(x)$ for $x\in[0,\frac{1}{2}]$. 
Note that \eqref{eq:gy:fracLaplace of u} holds also at $x=0,$ since $r-2s>0$ and both sides of the equation are equal to zero.
We thus obtain, by summing up 
$(-\Delta)^sv+(-\Delta)^s\varphi=(-\Delta)^su$, the following expression:
\begin{equation}\label{gy:eq:Deltas u with varphi}
   (-\Delta)^su(x)=c_1 x^{\beta}+c_2 x^{2s+\beta}+c_sG(x)+c_sH(x)\quad
   \text{for all }x\in {\textstyle[0,\frac{1}{2}]},
\end{equation}
where $G\in C^1([0,\frac{1}{2}])$, $G(0)=0$, and $H$ is here given by
\begin{equation}\label{gy:eq:H for u v varphi}
  H(x):=\int_{\R\setminus(-\frac{1}{2},\frac{1}{2})}dy\, 
  \frac{u(x)-u(x+y)}{|y|^{1+2s}}.
\end{equation}
We claim that the function $H$ is Lipschitz in $\R$. This easily follows from the fact that $u$ is itself Lipschitz (since $2s+\beta>1$). Suppose that $u$ has Lipschitz constant 
$M$, and thus 
\begin{equation}
|u(x)-u(x+y)-(u(z)-u(z+y))|\leq|u(x)-u(z)|+|u(x+y)-u(z+y)|
\leq 2M|x-z|.
\end{equation}
Combining this with the fact that $y\mapsto y^{-(1+2s)}$ is integrable over $(\frac{1}{2},+\infty)$ shows that $H$ is indeed Lipschitz in $\R$.

We use now that $u$ is invertible in $(-1,1)$ to define 
$f:[0,u(\frac{1}{2})]\to \R$  by
$$ 
   f(t):=c_1(u^{-1}(t))^{\beta}+c_2(u^{-1}(t))^{2s+\beta}+c_sG(u^{-1}(t))+c_sH(u^{-1}(t)).
$$
By \eqref{gy:eq:Deltas u with varphi}, $f$ and $u$ satisfy $(-\Delta)^su(x)=f(u(x))$ for all $x\in [0,\frac{1}{2}].$ 

Let us now show that $f\in C^{\beta}([0,u(\frac{1}{2})]).$  By the properties of $G$ and $H$ and the definition of $f$, it is enough to prove that $u^{-1}$ is Lipschitz in 
$[0,u(\frac{1}{2})]$. Let $0\leq p<q\leq u(\frac{1}{2})$ and set $p=u(x')$ and $q=u(y')$ for some $0\leq x'<y'\leq\frac{1}{2}$. Then, abbreviating $t'=x'/y'\leq 1$, we obtain
\begin{align}
  \frac{|u^{-1}(p)-u^{-1}(q)|}{|p-q|}=
  \frac{y'-x'}{u(y')-u(x')}=
  \frac{1-t'}{\left( (y')^{2s+\beta-1}+1 \right) -t' \left( (x')^{2s+\beta-1}+1 \right)}\leq 1.
\end{align}
The last inequality holds true for every $t'\in [0,1]$, not only for $t'=x'/y'$. 

To conclude the proof we now take the odd extension of $f$ to $[-u(\frac{1}{2}),0]=[u(-\frac{1}{2}),0]$, we proceed exactly as at the end of the proof of \Cref{gy:lemma:v more gener odd}, and we extend $f$ in an arbitrary $C^{\beta}$ way from 
$[-u(\frac{1}{2}),u(\frac{1}{2})]$ to the whole real line.
\end{proof}

We can now prove \Cref{gy:thm:intro:examples reg} $(ii)$. We only give a brief outline of the proof, as the construction is almost identical to the one made in the proof of \Cref{gy:example_8periodic:lemma}, but it is simpler. Notably, the fact that $H$ is Lipschitz in $\R$ is trivial (for a new periodized version of \eqref{gy:eq:u for ii nonper}) and follows exactly as in the proof of \Cref{prop:gy:counterexample reg ii} of the present section. Therefore, we will neither require property 
$(e)$ from \Cref{section:periodic counterexample i} nor a version of \Cref{lemma:gy:H property u periodic}.

\begin{proof}[Proof of  \Cref{gy:thm:intro:examples reg} $(ii)$]
We shall chose $u$ as an $8$-periodic function which coincides with the function defined in \eqref{gy:eq:u for ii nonper} in $[-1,1],$ is smooth in $\R\setminus 4\Z,$ increasing in $[0,2],$ odd with respect to $x=0$ and even with respect to $x=2,$ and is a quadratic function (in the sense of $(f)$ below \Cref{gy:example_8periodic:lemma}) in $(2-\delta,2+\delta)$ for some $\delta>0$.  

We first split the function $u$ into a sum $u=v+\varphi$ as in \eqref{gy:eq:split u in v plus varphi}, but where now $v$ and $\varphi$ are also $8$-periodic (this step is independent of how we define $v$ and $\varphi$ in $(1,2)$). Then, using \Cref{lemma:gy:fracLaplace of u} on $v$, we get that $u$ satisfies \eqref{gy:eq:Deltas u with varphi} and \eqref{gy:eq:H for u v varphi}. Proceeding as in the proof of \Cref{prop:gy:counterexample reg ii} we obtain a function $f\in C^{\beta}([0,u(\frac{1}{2})])$ which satisfies 
$(-\Delta)^su=f(u)$ in $[0,\frac{1}{2}].$ We can then conclude arguing as in the proof of  \Cref{gy:example_8periodic:lemma}, that is, extending $f$ to 
$(u(\frac{1}{2}),u(2)]$ by
$f(t):=((-\Delta)^su)(u^{-1}(t))$ and then taking the odd extension of $f$.
\end{proof}


\appendix

\section{An example of H\"older regularity for $u''=f(u)$}
\label{section:appendix opt reg local}

We give here a simple example regarding the optimal regularity of bounded solutions to the equation $u''=f(u)$ in some open interval $I$ with $f\in C^{\beta}(\R)$ for some 
$0<\beta<1$. We shall show that in general $u$ cannot be more regular than $C^{2+\beta}(I).$ Observe that taking the power function $u(x)=x^r$ to some power $r=r(\beta)$ will not give an appropriate example, as then $u$ turns out to be in $C^r(I)$ for some $r>2+\beta$. The following example that we will give inspired the construction in the proof of \Cref{gy:thm:intro:examples reg} $(ii)$ based on \eqref{gy:eq:u for ii nonper}. Moreover, we have not found in the literature any example of this sort and of this simplicity. 

Given $0<\beta<1$, define $u,f:\R\to\R$ by
$$
  u(x):=\left\{
  \begin{array}{rl}
       x^{2+\beta}+x &\text{if }x\geq 0,
       \\
       x & \text{if }x< 0,
  \end{array}\right.
  \qquad
    f(t):=\left\{
  \begin{array}{rl}
       (2+\beta)(1+\beta)(u^{-1}(t))^{\beta} &\text{if }t\geq 0,
       \\
       0 & \text{if }t < 0.
  \end{array}\right.
$$
Here $u^{-1}$ denotes the inverse of $u,$ which is obviously well defined in the whole real line. It is immediate to check that $u''=f(u)$ in $\R.$ It is also clear that $u\notin C^{2+\beta+\epsilon}([-\delta,\delta])$ for any $\epsilon,\delta>0.$ 
It only remains to show that $f\in C^{\beta}([-M,M])$ for every $M>0.$ Let $0\leq p<q\leq M$ and, without loss of generality, $p=u(x)$ and $q=u(y)$ with $0\leq x<y$. Setting $t=x/y<1$,  we obtain that 
\begin{align}
  \frac{|f(p)-f(q)|}{|p-q|^{\beta}}=
  (2+\beta)(1+\beta)\frac{1-t^{\beta}}{\left((y^{1+\beta}+1)-t(x^{1+\beta}+1)\right)^{\beta}}
  \leq (2+\beta)(1+\beta),
\end{align} 
as desired.
In the last inequality we have used that, for every $0\leq x<y$, the function $$t\in (0,1)\mapsto
1-t^{\beta}-((y^{1+\beta}+1)-t(x^{1+\beta}+1))^{\beta}$$ is negative due to convexity and the fact that it is negative at $t=0$ and at $t=1.$



\end{document}